\documentclass[12pt,draftcls,onecolumn]{IEEEtran}

\usepackage{epsfig} 
\usepackage[cmex10]{amsmath}
\usepackage{amssymb}
\usepackage[boxed]{algorithm2e}
\usepackage{cite}
\usepackage{graphicx}

\newtheorem{theorem}{Theorem}
\newtheorem{lemma}{Lemma}
\newtheorem{assumption}{Assumption}

\newcommand{\optvar}{ {(\text{opt})} }
\newcommand{\dvar}{ {\lambda} }
\newcommand{\smooth}{ {\text{s}} }
\newcommand{\poly}{ {\text{p}} }
\newcommand{\tr}{ \top }

\title{Time-Average Optimization with Non-Convex Decision Set and Its Convergence}

\author{Sucha Supittayapornpong, Longbo Huang, Michael J. Neely
\thanks{This material is supported in part by one or more of: the NSF Career grant CCF-0747525, the Network Science Collaborative Technology Alliance sponsored by the U.S. Army Research Laboratory W911NF-09-2-0053, the National Basic Research Program of China Grant 2011CBA00300, 2011CBA00301, the National Natural Science Foundation of China Grant 61033001, 61361136003, 61303195, the China youth 1000-talent grant.}
\thanks{S. Supittayapornpong and M. J. Neely are with Electrical Engineering Department, University of Southern California, 3740 McClintock Ave., Los Angeles, CA, USA 90089-2565, Tel: +1-213-740-4685, Fax: +1-213-740-8729
        {\tt\small supittay@usc.edu, mjneely@usc.edu}}%
\thanks{L. Huang is with Institute for Interdisciplinary Information Sciences, Tsinghua University, Beijing, China, 100084, Tel: +8610-62781693, Fax: +8610-62797331-2000
        {\tt\small longbohuang@tsinghua.edu.cn}}%
\thanks{The corresponding author is Mr. Sucha Supittayapornpong.}
}

\newcommand{\nth}{ {\text{\tiny th}} }

\newcommand{\set}[1]{ {\mathcal{#1}} }

\newcommand{\abs}[1]{ {\left| {#1} \right|} }

\newcommand{\norm}[1]{ {\left\lVert {#1} \right\rVert} }
\newcommand{\RealSet}{ {\mathbb{R}} }

\newcommand{\defeq}{ {\triangleq} }

\newcommand{\arginf}{\operatornamewithlimits{arginf}}

\newcommand{\maximize}{ {\text{Maximize}} }
\newcommand{\minimize}{ {\text{Minimize}} }
\newcommand{\subjectto}{ {\text{Subject to}} }

\newcommand{\prts}[1]{ {\left[ {#1} \right]} }
\newcommand{\prtr}[1]{ {\left( {#1} \right)} }
\newcommand{\prtc}[1]{ {\left\{ {#1} \right\}} }

\newcommand{\prtbs}[1]{ {\biggl[ {#1} \biggr]} }

\newcommand{\minvar}{ {(\text{min})} }

\begin{document}
\maketitle


\begin{abstract}
This paper considers \textit{time-average optimization}, where a decision vector is chosen every time step within a (possibly non-convex) set, and the goal is to minimize a convex function of the time averages subject to convex constraints on these averages.  Such problems have applications in networking, multi-agent systems, and operations research, where decisions are constrained to a discrete set and the decision average can represent average bit rates or average agent actions.  This time-average optimization extends traditional convex formulations to allow a non-convex decision set.  This class of problems can be solved by Lyapunov optimization.  A simple drift-based algorithm, related to a classical dual subgradient algorithm, converges to an $\epsilon$-optimal solution within $O(1/\epsilon^2)$ time steps.  Further, the algorithm is shown to have a transient phase and a steady state phase which can be exploited to improve convergence rates to $O(1/\epsilon)$ and $O(1/{\epsilon^{1.5}})$ when vectors of Lagrange multipliers satisfy locally-polyhedral and locally-smooth assumptions respectively.  Practically, this improved convergence suggests that decisions should be implemented after the transient period.
\end{abstract}

\section{INTRODUCTION}
\label{sec:intro}
Convex optimization is often used to optimally control communication networks (see \cite{Chiang:Layered} and references therein) and distributed multi-agent systems \cite{Nedic:multi_agent}.  This framework utilizes both convexity properties of an objective function and a feasible decision set.  However, various systems have inherent discrete (and hence non-convex) decision sets.  For example, a wireless system might constrain transmission rates to a finite set corresponding to a fixed set of coding options.  Further, distributed agents might only have finite options of decisions.  This discreteness restrains the application of convex optimization.

Let $I$ and $J$ be positive integers. This paper considers a class of problems called \textit{time-average optimization} where decision vectors $x(t) = (x_1(t), \dotsc, x_I(t))$ are chosen sequentially over time slots $t \in \{0, 1, 2, \dotsc\}$ from a decision set $\set{X}$, which is a closed and bounded subset of $\mathbb{R}^I$ (possibly non-convex and discrete), and its average $\bar{x} = \lim_{T \rightarrow \infty} \frac{1}{T} \sum_{t = 0}^{T-1} x(t)$ solves the following problem:
\begin{align}
  \label{eq:ori_formulation}
  \minimize \quad & f \prtr{ \bar{x} } \\
  \subjectto \quad & g_j \prtr{ \bar{x} } \leq 0 && j \in \{ 1, \dotsc, J \} \notag\\
  {} & x(t) \in \set{X} && t \in \{ 0, 1, 2, \dotsc \}, \notag
\end{align}
where $f: \overline{\set{X}} \rightarrow \RealSet$ and $g_j: \overline{\set{X}} \rightarrow \RealSet$ are convex functions and $\overline{\set{X}}$ is the convex hull of $\set{X}$.

This time-average optimization reflects a scenario where an objective is in the time-average sense.  For example, network users are interested in average bit rates or throughput, and distributed agents are concerned with average actions.  The formulation can be considered as a fine granularity version of a one-shot average formulation, where an average decision is chosen, and can be used to extend several convex optimization problems in literature, see for example \cite{Chiang:Layered} and references therein, to have non-convex decision sets.  



Formulation \eqref{eq:ori_formulation} has an optimal solution which can be converted (by averaging) to the following convex optimization problem:
\begin{align}
  \label{eq:convex_formulation}
  \minimize \quad & f \prtr{ x } \\
  \subjectto \quad & g_j \prtr{ x } \leq 0 && j \in \{ 1, \dotsc, J \} \notag\\
  {} & x \in \overline{\set{X}}. \notag
\end{align}
Note that an optimal solution to formulation \eqref{eq:convex_formulation} may not be in the non-convex decision set $\set{X}$. Nevertheless, problems \eqref{eq:ori_formulation} and \eqref{eq:convex_formulation} have the same optimal value.  In addition, directly applying a primal-average technique on a non-convex formation \eqref{eq:non_convex_formulation}, where the convex hull in \eqref{eq:convex_formulation} is removed, may lead to an local optimal solution with respect to the time-average problem \eqref{eq:ori_formulation}.  For example, when $\set{X} = \prtc{0,1}, J=1, f(x) = ( x - 2/3 )^2, g_1(x) = 2/3 - x$, a primal average solution of the technique in \cite{Nedic:approximate_primal} is $1$, while a solution to problem \eqref{eq:ori_formulation} is $\bar{x} = 2/3$.

\begin{align}
  \label{eq:non_convex_formulation}
  \minimize \quad & f \prtr{ x } \\
  \subjectto \quad & g_j \prtr{ x } \leq 0 && j \in \{ 1, \dotsc, J \} \notag\\
  {} & x \in \set{X}. \notag
\end{align}

Although there have been several techniques utilizing time-average solutions \cite{Nesterov:dual_averaging,Nedic:approximate_primal, neely-dist-comp}, those works are limited to convex formulations.  In fact, this work can be considered as a generalization of \cite{Nedic:approximate_primal, neely-dist-comp} as decisions are allowed to be chosen from a non-convex set.  A non-convex optimization problem is considered in \cite{Minghui:nonconvex}, where an approximate problem is solved with the assumption of a unique vector of Lagrange multipliers.  In comparison, when $f(x)$ and $g_j(x)$'s are Lipschitz continuous, the algorithm proposed in this paper solves problem \ref{eq:ori_formulation} without the uniqueness assumption.   This paper is inspired by the Lyapunov optimization technique \cite{Neely:SNO} which solves stochastic and time-average optimization problems, including problems such as \eqref{eq:ori_formulation}.  This paper removes the stochastic characteristic and focuses on the connection between the technique and a general convex optimization.  This allows a convergence time analysis of a \emph{drift-plus-penalty} algorithm that solves problem \eqref{eq:ori_formulation}.  Importantly, this paper shows that faster convergence can be achieved by starting time averages after a suitable transient period.

Another area of literature focuses on convergence time of first-order algorithms to an $\epsilon$-optimal solution to a convex problem, including problem \eqref{eq:convex_formulation}.  For \emph{unconstrained optimization} without strong convexity of the objective function, the accelerated method (with Lipschitz continuous gradients) has $O(1/\sqrt{\epsilon})$ convergence time \cite{Nesterov:O,Tseng:Fast_gradient}, while gradient and subgradient methods take $O(1/\epsilon)$ and $O(1/\epsilon^2)$ respectively \cite{Boyd:CO,Nedic:approximate_primal}.  Two $O(1/\epsilon)$ first-order methods for \emph{constrained optimization} are developed in \cite{Beck:FastNUM,Ermin:AsynADMM}, but the results rely on special convex formulations.  A second-order method for constrained optimization \cite{Xia:Newton} has a fast convergence rate but relies on special a convex formulation.  All of these results rely on convexity assumptions that do not hold in formulation \eqref{eq:ori_formulation}.

This paper develops an algorithm for the formulation \eqref{eq:ori_formulation} and analyzes its convergence time.  The algorithm is shown to have $O(1/\epsilon^2)$ convergence time with a mild Slater condition.  However, inspired by results in \cite{Longbo:Delay_reduction}, under a uniqueness assumption on Lagrange multipliers the algorithm is shown to enter two phases: a \emph{transient phase}  and a \emph{steady state phase}.  Convergence time can be significantly improved by starting the time averages after the transient phase.  Specifically, when a dual function satisfies a \textit{locally-polyhedral} assumption, the modified algorithm has $O(1/{\epsilon})$ convergence time (including the time spent in the transient phase), which equals the best known convergence time for constrained convex optimization via first-order methods.  On the other hand, when the dual function satisfies a \textit{locally-smooth} assumption, the algorithm has $O(1/{\epsilon^{1.5}})$ convergence time.  Furthermore, simulations show that these fast convergence times are robust even without the uniqueness assumption.  An application of these improved convergence times can be effective implementation of decisions where decisions are implemented online after offline calculation during a transient period.

The contributions of this paper are summarized below.
\begin{enumerate}
\item We establish the connection between Lyapunov optimization and a dual subgradient algorithm for a problem with a non-convex decision set, which requires additional problem transformation.
\item We generalize the modeling of a one-shot convex optimization \eqref{eq:convex_formulation}, extensively used in \cite{Chiang:Layered}, to the time-average formulation \eqref{eq:ori_formulation} that allows a non-convex decision set, while optimality and complexity are preserved.
\item We investigate transient and steady-state behaviors of the algorithm solving the time-average problem \eqref{eq:ori_formulation}.  Then, we exploit the behaviors to obtain sequences of decisions that achieve $O(\epsilon)$-optimal solutions within $O(1/\epsilon)$ and $O(1/{\epsilon^{1.5}})$ iterations under locally-polyhedral and locally-smooth assumptions instead of the standard $O(1/{\epsilon^{2}})$ iterations in \cite{Nedic:approximate_primal, neely-dist-comp}.
\end{enumerate}

The paper is organized as follows.  Section \ref{sec:DPP} constructs an algorithm to solve the time-average problem.  The general $O(1/\epsilon^2)$ convergence time is proven in Section \ref{section:convergence}.  Section \ref{sec:fast_rates} explores faster convergence times of $O(1/\epsilon)$ and $O(1/{\epsilon^{1.5}}$) under the unique Lagrange multiplier assumption.  Example problems are given in Section \ref{sec:example}, including cases when the uniqueness condition fails.  Section \ref{sec:conclusion} concludes the paper.

\section{TIME-AVERAGE OPTIMIZATION}
\label{sec:DPP}

In order to solve problem \eqref{eq:ori_formulation}, an embedded problem with a similar solution is formulated with the following assumptions.

\subsection{The extended set $\set{Y}$} 
\label{sec:sety}

Let $\set{Y}$ be a closed, bounded, and convex subset of $\mathbb{R}^I$ that contains $\overline{\set{X}}$.
Assume the functions $f(x)$, $g_j(x)$ for $j \in \{1, \ldots, J\}$ extend as real-valued convex functions over $x \in \set{Y}$.  
The set $\set{Y}$ can be defined as $\overline{\set{X}}$ itself.  
However, choosing $\set{Y}$ as a larger set helps to ensure a Slater condition is satisfied (defined below).   
Further, choosing $\set{Y}$ to have a simple structure helps to simplify the resulting optimization.  
For example, set $\set{Y}$ might be chosen as a closed and bounded hyper-rectangle that contains $\overline{\set{X}}$ in its interior.

\subsection{Lipschitz continuity and Slater condition} 

In addition to assuming that $f(x)$ and $g_j(x)$ are convex over $x \in \set{Y}$, assume they are \emph{Lipschitz continuous}, so there is a constant $M>0$ such that for all $x, y \in \set{Y}$:
\begin{align}
  \abs{ f(x) - f(y) } & \leq M \norm{x - y}  \label{eq:lip1} \\
  \abs{ g_j(x) - g_j(y) } &\leq  M \norm{x - y} \label{eq:lip2} 
\end{align}
where $\norm{x} = \sqrt{x_1^2 + \cdots + x_I^2}$ is the Euclidean norm.

Further, assume that there exists a vector $\hat{x} \in \overline{\set{X}}$ that satisfies $g_j(\hat{x}) < 0$ for all $j \in \{1, \ldots, J\}$, and is such that $\hat{x}$ is in the interior of set $\set{Y}$.  
This is a \emph{Slater condition} that, among other things, ensures the constraints are feasible for the problem of interest.

\subsection{Relation to dual subgradient algorithm}

Problem \eqref{eq:ori_formulation} can be solved by the Lyapunov optimization technique \cite{Neely:SNO}.  It has been known that the drift-plus-penalty algorithm in the Lyapunov optimization is identical to a classic dual subgradient method \cite{Bertsekas:Convex,Nedic:approximate_primal} that solves problem \eqref{eq:emb_formulation}, with the exception that it takes a time average of primal values.  
\begin{align}
  \label{eq:emb_formulation}
  \minimize \quad & f(y) \\
  \subjectto \quad & g_j(y) \leq 0 && j \in \{ 1, \dotsc, J \} \notag \\
  {} & x_i = y_i && i \in \{ 1, \dotsc, I \} \notag \\
  {} & x \in \overline{\set{X}}, \quad y \in \set{Y}. \notag
\end{align}
This was noted in \cite{Neely:Power,Longbo:Delay_reduction} for related problems.  Problem \eqref{eq:emb_formulation} is called the \emph{embedded formulation} of the time-average problem \eqref{eq:ori_formulation} and is convex.  It is not difficult to show that the above problem has an optimal value $f^\optvar$ that is the same as that of problems \eqref{eq:ori_formulation} and \eqref{eq:convex_formulation}.  Compared to a formulation in \cite{Nedic:approximate_primal}, problem \eqref{eq:emb_formulation} contains additional equality constraints and $\overline{\set{X}}$ derived from the original decision set.  This makes further analysis and algorithm slightly different from \cite{Nedic:approximate_primal}, whose results cannot be applied directly.

Now consider the dual of embedded formulation \eqref{eq:emb_formulation}.  
Let vectors $w$ and $z$ be dual variables of the first and second constraints in problem \eqref{eq:emb_formulation}, where the feasible set of $(w, z)$ is denoted by $\Pi = \RealSet_+^J \times \RealSet^I$.  Let $g(y) = (g_1(y), \dotsc, g_J(y))$ denote a $J$-dimensional column vector of functions $g_j(y)$.
A Lagrangian has the following expression:
\begin{equation*}
  \Lambda(x, y, w, z) = f(y) + w^\tr g(y) + z^\tr (x-y).
\end{equation*}
Define:
\begin{align*}
  x^\ast(z) & = \arginf_{x \in \overline{\set{X}}} z^\tr x \quad (\text{with}~x^\ast(z) \in \set{X})\\
  y^\ast(w, z) & = \arginf_{y \in \set{Y}} [f(y) + w^\tr g(y) - z^\tr y].
\end{align*}
Notice that $x^\ast(z)$ may have multiple candidates including extreme point solutions, since $z^\tr x$ is a linear function.
We restrict $x^\ast(z)$ to any of these extreme solutions, which implies $x^\ast(z) \in \set{X}$.  
Then the dual function is defined as
\begin{align}
  \label{eq:dual_function}
  d(w, z) &= \inf_{x \in \overline{\set{X}}, y \in \set{Y}} \Lambda(x, y, w, z) \\
  & \hspace{-0.5em}\quad = f(y^\ast(w,z)) + w^\tr g(y^\ast(w, z)) + z^\tr [ x^\ast(z) - y^\ast(w, z) ]. \notag
\end{align}
A pair of subgradients \cite{Bertsekas:Convex} with respect to $w$ and $z$ is:
\begin{equation*}
  \partial_w d(w, z) = g(y^\ast(w, z)), \quad  \partial_z d(w, z)  = x^\ast(z) - y^\ast(w, z).
\end{equation*}

Finally, the dual formulation of embedded problem \eqref{eq:emb_formulation} is
\begin{align}
  \label{eq:dual_formulation}
  \maximize \quad & d(w, z) \\
  \subjectto \quad & (w, z) \in \Pi. \notag
\end{align}
Let the optimal value of problem \eqref{eq:dual_formulation} be $d^\ast$.  
Since problem \eqref{eq:emb_formulation} is convex, the duality gap is zero, and $d^\ast = f^\optvar$.  
Problem \eqref{eq:dual_formulation} can be treated by a dual subgradient method \cite{Bertsekas:Convex} with a fixed stepsize $1/V$ and the restriction on $x(t) \in \set{X}$, where $V > 0$ is a parameter.  This leads to Algorithm \ref{alg:Subgrad} summarized in the figure below, called the \emph{dual subgradient algorithm}.  Note that the algorithm is different from the one in \cite{Nedic:approximate_primal} due to the equality constraints and the restriction on $x(t)$.
\begin{algorithm}
  \DontPrintSemicolon
  Initialize $w(0)$ and $z(0)$.\;
  \For{ $t = 0, 1, 2, \dotsc $ } {
    $x(t) = \arginf_{x \in \overline{\set{X}}} z(t)^\tr x$ \quad(with $x(t) \in \set{X}$)\;
    $y(t) = \arginf_{y \in \set{Y}} [f(y) + w(t)^\tr g(y) - z(t)^\tr y]$\;
    $w(t+1) = \prts{ w(t) + \frac{1}{V} g(y(t)) }_+$\;
    $z(t+1) = z(t) + \frac{1}{V} [ x(t) - y(t) ]$\;
  }
  \caption{Dual subgradient algorithm with restriction}
  \label{alg:Subgrad}
\end{algorithm}


Traditionally, the dual subgradient algorithm of \cite{Bertsekas:Convex} is intended to produce primal vector estimates that converge to a desired result. 
However, this requires additional assumptions.  
Indeed, for our problem, the primal vectors $x(t)$ and $y(t)$ \emph{do not converge} to anything near a solution in many cases, such as when the $f(x)$ and $g_j(x)$ functions are linear or piecewise linear.
However, Algorithm \ref{alg:Subgrad} ensures that the \emph{time averages} of $x(t)$ and $y(t)$ converge as desired.  

We use the notation $w(t)$ and $z(t)$ from Algorithm \ref{alg:Subgrad}, with the update rule for $w(t+1)$ and $z(t+1)$ given there:
\begin{align}
  w(t+1) & = \prts{ w(t) + \frac{1}{V} g(y(t)) }_+ \label{eq:w_update}\\
  z(t+1) & = z(t) + \frac{1}{V} [ x(t) - y(t) ]. \label{eq:z_update}
\end{align}
For ease of notation, define $\dvar(t) \defeq (w(t), z(t))$ as a concatenation of these vectors.  Let $C$ be some positive constant such that $\norm{ g(y) }^2 \leq C$ and $\norm{ x - y }^2 \leq C$ for any $x \in \overline{\set{X}}$ and any $y \in \set{Y}$, since $\overline{\set{X}}$ is closed and bounded.  We first provide some useful properties.  It holds that
\begin{equation}
  \label{eq:onestep_bound}
  \norm{ \dvar(t+1) - \dvar(t) } \leq \sqrt{2C}/V \quad \text{for all}~ t,
\end{equation}
since
\begin{align}
  \norm{ \dvar(t+1) - \dvar(t) }^2
  & = \norm{ w(t+1) - w(t) }^2 + \norm{ z(t+1) - z(t) }^2 \notag\\
  & \leq \frac{1}{V^2} \norm{ g(y(t)) }^2 + \frac{1}{V^2} \norm{ x(t) - y(t) }^2 \label{eq:dvardrift1}\\
  & \leq 2C/V^2 \label{eq:dvardrift2}
\end{align}
where \eqref{eq:dvardrift1} follows from equations \eqref{eq:w_update}--\eqref{eq:z_update}, and \eqref{eq:dvardrift2} follows from the definition of $C$. Further,
\begin{align*}
  \norm{\dvar(t+1)}^2 - \norm{\dvar(t)}^2
  & = \norm{w(t+1)}^2 + \norm{z(t+1)}^2 - \norm{w(t)}^2 - \norm{z(t)}^2 \notag\\
  & \leq \frac{2C}{V^2} + \frac{2}{V}w(t)^\tr g(y(t)) + \frac{2}{V}z(t) \prts{x(t) - y(t)},
\end{align*}
where the last inequality uses the result of expanding the square norms of \eqref{eq:w_update} and \eqref{eq:z_update}.  Since Algorithm \ref{alg:Subgrad} chooses $x(t)$, $y(t)$ to minimize $d(\dvar(t)) = d(w(t),z(t))$ in \eqref{eq:dual_function}, the above bound and \eqref{eq:dual_function} imply that
\begin{align}
  d(\dvar(t)) &= f(y(t)) + w(t)^\tr g(y(t)) + z(t)^\tr\prts{x(t) - y(t)} \notag\\
  & \geq f(y(t)) + \frac{V}{2}\prts{ \norm{\dvar(t+1)}^2 - \norm{\dvar(t)}^2 } - \frac{C}{V}. \label{eq:modified_DPP}
\end{align}

From convex analysis, the dual function $d(\dvar)$, defined in \eqref{eq:dual_function}, has the following properties \cite{Bertsekas:Convex}:

\begin{itemize}
\item $d(\dvar)\leq f^\optvar$ for all $\dvar \in \Pi$.

\item If the Slater condition holds, then there are real numbers $F>0$,  $\eta > 0$ such that: 
\begin{equation*}
  d(\lambda) \leq F - \eta \norm{\lambda} \quad \text{for all}~ \dvar \in \Pi.
\end{equation*}

\item If the Slater condition holds, then there is an optimal value $\dvar^\ast \in \Pi$, called a \emph{Lagrange multiplier vector} \cite{Bertsekas:Convex}, that maximizes $d(\dvar)$.  Specifically, $d(\dvar^\ast) = f^\optvar$.
\end{itemize}



The first two properties can be substituted into the inequality \eqref{eq:modified_DPP} to ensure that, under Algorithm \ref{alg:Subgrad}, the following inequalities hold for all time slots $t \in \{0, 1, 2, \dotsc \}$:
\begin{align}
  \frac{V}{2}\left[\norm{\lambda(t+1)}^2 - \norm{\lambda(t)}^2\right] + f(y(t)) 
  & \leq \frac{C}{V} + f^\optvar \label{eq:DPP_opt} \\
  \frac{V}{2}\left[\norm{\lambda(t+1)}^2 - \norm{\lambda(t)}^2\right] + f(y(t))  
  & \leq \frac{C}{V} + F - \eta\norm{\lambda(t)} \label{eq:DPP_slater}
\end{align}

\section{GENERAL CONVERGENCE RESULT} \label{section:convergence} 

Define the average of variables $\{ a(t) \}_{t=0}^{T-1}$ as
\begin{equation*}
  \bar{a}(T) \defeq \frac{1}{T} \sum_{t=0}^{T-1} a(t) \quad\quad \text{for } T \in \prtc{1, 2, \dotsc }.
\end{equation*}

\begin{theorem}
  \label{thm:gen_conv}
  Let $\{ x(t), w(t), z(t) \}_{t=0}^\infty$ be a sequence generated by Algorithm \ref{alg:Subgrad}.  For $T > 0$, we have
  \begin{equation}
    \label{eq:avg_bound}
    f(\bar{x}(T)) - f^\optvar \leq \frac{V}{2T} \prts{ \norm{\dvar(0)}^2 - \norm{\dvar(T)}^2 } + \frac{C}{V} + \frac{VM}{T} \norm{ z(T) - z(0) } 
\end{equation}
\vspace{-1.5em}
  \begin{equation}
    \label{eq:avg_vio}
    g_j(\bar{x}(T)) \leq \frac{V}{T} \abs{ w_j(T) - w_j(0) } + \frac{V M}{T} \norm{ z(T) - z(0) } \\\quad j \in \{ 1, \dotsc, J \},
  \end{equation}
  where $M$ is the Lipschitz constant from \eqref{eq:lip1}--\eqref{eq:lip2}. 
\end{theorem}
\begin{IEEEproof}
For the first part, we have from the Lipschitz property \eqref{eq:lip1}: 
\begin{equation}
  \label{genbound_1}
   f(\bar{x}(T)) - f^\optvar  \leq [ f(\bar{y}(T)) - f^\optvar ] +  M \norm{ \bar{y}(T) - \bar{x}(T) }.
\end{equation}
We first upper bound $f(\bar{y}(T)) - f^\optvar$ on the right-hand side of \eqref{genbound_1}.  Let $\{ x(t), y(t), w(t), z(t) \}_{t=0}^\infty$ be a sequence generated by Algorithm \ref{alg:Subgrad}.  Relation \eqref{eq:DPP_opt} can be rewritten as
\begin{equation*}
  f(y(t)) - f^\optvar \leq \frac{C}{V} + \frac{V}{2} \prts{ \norm{\dvar(t)}^2 - \norm{\dvar(t+1)}^2 }.
\end{equation*}
Summing from $t = 0, \dotsc, T-1$ and dividing by $T$ give:
\begin{equation*}
  \frac{1}{T} \sum_{t=0}^{T-1} f(y(t)) - f^\optvar \leq \frac{C}{V} + \frac{V}{2T} \prts{ \norm{\dvar(0)}^2 - \norm{\dvar(T)}^2 }.
\end{equation*}
Using Jensen's inequality and the convexity of $f(\cdot)$ give:
\begin{equation}
  \label{eq:opt_bound}
  f(\bar{y}(T)) - f^\optvar \leq \frac{V}{2T} \prts{ \norm{\dvar(0)}^2 - \norm{\dvar(T)}^2 } + \frac{C}{V}.
\end{equation}
For $\norm{ \bar{y}(T) - \bar{x}(T) }$ in \eqref{genbound_1}, we consider the update equation of $z(t)$ in \eqref{eq:z_update}.  Summing from $t = 0, \dotsc, T-1$ yields $z_i(T) - z_i(0) = \frac{1}{V} \sum_{t=0}^{T-1} \prts{ x_i(t) - y_i(t) }$ for every $i$.  Rearranging and dividing by $T$ gives:
\begin{equation}
  \bar{x}_i(T) - \bar{y}_i(T) = \frac{V}{T} \prts{ z_i(T) - z_i(0) } \quad i \in \{ 1, \dotsc, I \}. \label{eq:vio_xy}
\end{equation}
Substituting \eqref{eq:opt_bound} and \eqref{eq:vio_xy} into \eqref{genbound_1} proves \eqref{eq:avg_bound}.

For the second part, we have from \eqref{eq:lip2}: 
\begin{equation}
  \label{genbound_2}
  g_j(\bar{x}(T)) \leq g_j(\bar{y}(T)) + M \norm{ \bar{y}(T) - \bar{x}(T) }.
\end{equation}
We first bound $g_j(\bar{y}(T))$.  The update equation of $w(t)$ in \eqref{eq:w_update} implies, for every $j$, that
\begin{equation*}
  w_j (t+1) = [ w_j (t) + \frac{1}{V} g_j (y(t)) ]_+ \geq w_j(t) + \frac{1}{V} g_j (y(t)),
\end{equation*}
and $w_j(t+1) - w_j (t) \geq \frac{1}{V} g_j (y(t))$.  Summing from $t = 0, \dotsc, T-1$, we have $w_j(T) - w_j(0) \geq \frac{1}{V} \sum_{t=0}^{T-1} g_j (y(t))$.  
Dividing by $T$ and using Jensen's inequality and convexity of $g_j( \cdot )$ gives
\begin{equation*}
  \frac{1}{T} [ w_j (T) - w_j (0) ] \geq \frac{1}{VT} \sum_{t=0}^{T-1} g_j (y(t)) \geq \frac{1}{V} g_j ( \bar{y}(T) ).
\end{equation*}
This shows that
\begin{equation}
    g_j( \bar{y}(T) ) \leq \frac{V}{T} \abs{ w_j(T) - w_j(0) } \quad j \in \{ 1, \dotsc, J \}. \label{eq:vio_g}
\end{equation}
Substituting \eqref{eq:vio_g} and \eqref{eq:vio_xy} into \eqref{genbound_2} proves \eqref{eq:avg_vio}.
\end{IEEEproof}

Theorem \ref{thm:gen_conv} can be interpreted when $\norm{ \dvar }$ is bounded from above by some finite constant as that the deviation from optimality \eqref{eq:avg_bound} is bounded from above by $O(V/T + 1/V)$, and the constraint violation \eqref{eq:avg_vio} is bounded above by $O(V/T)$.  To have both bounds be within $O(\epsilon)$, we set $V = 1/\epsilon$ and $T = 1/\epsilon^2$.  Thus the convergence time of Algorithm \ref{alg:Subgrad} is $O(1/\epsilon^2)$.  The next lemma shows that such a constant exists when the Slater condition holds.

\begin{lemma}
  \label{lem:l2} 
  When $V \geq 1$, $w_j(0) = z_i(0) = 0$ for all $i$ and $j$, then under Algorithm \ref{alg:Subgrad}, the Slater condition implies there is a constant $D > 0$ (independent of $V$) such that 
\begin{equation*}
  \norm{ \lambda(t) } = \sqrt{ \sum_{j=1}^J w_j(t)^2 + \sum_{i=1}^I z_i(t)^2 } \leq D \quad \text{for all } t.
\end{equation*}
\end{lemma}
\begin{IEEEproof}
  From \eqref{eq:DPP_slater} and $V \geq 1$, if $\norm{\dvar(t)} \geq (C + F - f^\minvar)/\eta$ where $f^\minvar = \inf_{y \in \set{Y}} f(y)$, then we have
\begin{align*}
  \frac{V}{2} \prts{ \norm{\dvar(t+1)}^2 - \norm{\dvar(t)}^2 } 
  & \leq \frac{C}{V} + F - f(y(t)) - \eta\norm{\dvar(t)} \\
  & \leq 0
\end{align*}
This implies that:
\begin{equation*}
  \norm{ \dvar(t) } \leq (C + F - f^\minvar)/\eta + \norm{ \dvar(t+1) - \dvar(t) }.
\end{equation*}
To complete the proof, note that $\norm{\dvar(t+1) - \dvar(t)} \leq \sqrt{2C}/V$ from \eqref{eq:onestep_bound}.
Since $V \geq 1$, letting $D \defeq (C + F - f^\minvar)/\eta + \sqrt{2C}$ proves the lemma.
\end{IEEEproof} 

This section shows that Algorithm \ref{alg:Subgrad} generates a sequence of decisions that achieves $O(\epsilon)$-optimal solution within $O(1/\epsilon^2)$ iterations.  The next section shows that it is possible to generate an $O(\epsilon)$-optimal achieving sequence of decisions within $O(1/\epsilon)$ and $O(1/\epsilon^{1.5})$ by analyzing a \emph{transient phase} and a \emph{steady state phase} of Algorithm \ref{alg:Subgrad}.

\section{CONVERGENCE OF TRANSIENT AND STEADY STATE PHASES}
\label{sec:fast_rates}

With this idea, we analyze the convergence time in the case when the dual function satisfies a \emph{locally-polyhedral} assumption and the case when it satisfies a \emph{locally-smooth} assumption.  Both cases use the following mild assumption:

\begin{assumption}
  \label{ass:uniqueness}
  The dual formulation \eqref{eq:dual_formulation} has a unique Lagrange multiplier denoted by $\dvar^\ast \defeq (w^\ast, z^\ast)$.
\end{assumption}

This assumption is assumed throughout Section \ref{sec:fast_rates}, and replaces the Slater assumption (which is no longer needed).  
Note that this is a mild assumption when practical systems are considered, e.g., \cite{Longbo:Delay_reduction,Eryilmaz:QLB}.  
In addition, simulations in Section \ref{sec:example} suggest that the algorithm derived in this section still has desirable performance without this uniqueness assumption.

We first provide a general result that will be used later.

\begin{lemma}
  \label{lem:dual_move}
  Let $\{ \dvar(t) \}_{t=0}^\infty$ be a sequence generated by Algorithm \ref{alg:Subgrad}.
  The following relation holds:
  \begin{equation}
    \label{eq:dual_move}
    \norm{\dvar(t+1) - \dvar^\ast}^2 \leq \norm{ \dvar(t) - \dvar^\ast }^2 + \frac{2}{V} [d(\dvar(t)) - d(\dvar^\ast)]\\
    + \frac{2 C}{V^2}, \quad t \in \prtc{0, 1, 2 , \dotsc}.
  \end{equation}
\end{lemma}
\begin{IEEEproof}
Recall that $\dvar(t) = (w(t), z(t))$. Define $h(t) \defeq (g(y(t)), x(t) - y(t))$ as the concatenation vector of the constraint functions.
From the non-expansive property, we have that
\begin{align}
  \norm{\dvar(t+1) - \dvar^\ast}^2
  & = \norm{ \prtr{ \prts{ w(t) + \frac{1}{V} g(y(t)) }_+, z(t) + \frac{1}{V}\prts{ x(t) - y(t) }} - \dvar^\ast }^2 \notag\\
  & \leq \norm{ \prtr{ w(t) + \frac{1}{V} g(y(t)) , z(t) + \frac{1}{V}\prts{ x(t) - y(t) }} - \dvar^\ast }^2 \notag\\
  & = \norm{ \dvar(t) + \frac{1}{V} h(t) - \dvar^\ast }^2 \notag\\
  & = \norm{ \dvar(t) - \dvar^\ast }^2 + \frac{1}{V^2} \norm{h(t)}^2 + \frac{2}{V} [\dvar(t) - \dvar^\ast]^\tr h(t) \notag\\
  & \leq \norm{ \dvar(t) - \dvar^\ast }^2 + \frac{2 C}{V^2} + \frac{2}{V} [d(\dvar(t)) - d(\dvar^\ast)],
\end{align}
where the last inequality uses the definition of $C$ and the concavity of the dual function \eqref{eq:dual_function}, i.e, $d(\dvar_1) \leq d(\dvar_2) + \partial d(\dvar_2)^\tr [\dvar_1 - \dvar_2]$ for any $\dvar_1, \dvar_2 \in \Pi$, and $\partial d(\dvar(t)) = h(t)$.
\end{IEEEproof}

\subsection{Locally-Polyhedral Dual Function}
\label{sec:poly}

\begin{figure}
  \centering
  \includegraphics[scale=0.8]{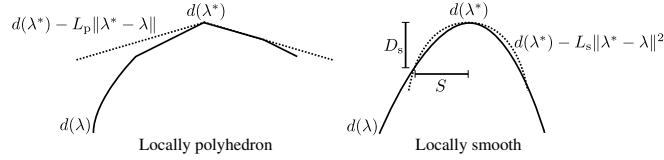}
  \caption{Illustration of locally-polyhedral and locally-smooth functions}
  \label{fig:ps_functions}
\end{figure}

Throughout Section \ref{sec:poly}, the dual function \eqref{eq:dual_function} is assumed to have a locally-polyhedral property, introduced in \cite{Longbo:Delay_reduction}, as stated in Assumption \ref{ass:poly}.  A dual function with this property is illustrated in Figure \ref{fig:ps_functions}.  The property holds when $f$ and $g_j$ for every $j$ are either linear or piece-wise linear.
\begin{assumption}
  \label{ass:poly}
  There exists an $L_\poly > 0$ such that the dual function \eqref{eq:dual_function} satisfies
\begin{equation}
  \label{eq:poly}
  d(\dvar^\ast) \geq d(\dvar) + L_\poly \norm{ \dvar - \dvar^\ast } \quad \text{for all}~\dvar \in \Pi
\end{equation}
where $\dvar^\ast$ is the unique Lagrange multiplier.
\end{assumption}

The ``p'' subscript in $L_\poly$ represents ``polyhedral.''  Furthermore, concavity of dual function \eqref{eq:dual_function} ensures that if this property holds locally about $\dvar^\ast$, it also holds globally for all $\dvar \in \Pi$ (see Figure \ref{fig:ps_functions}).

The behavior of the generated dual variables with dual function satisfying the locally-polyhedral assumption can be described as follows.  Define
\begin{equation*}
B_\poly(V) \defeq \max\prtc{ \frac{L_\poly}{2V}, \frac{2 C}{V L_\poly} }.
\end{equation*}
\begin{lemma}
  \label{lem:poly_drift}
  Under Assumptions \ref{ass:uniqueness} and \ref{ass:poly}, whenever $\norm{ \dvar(t) - \dvar^\ast } \geq B_\poly(V)$, it follows that
\begin{equation}
  \label{eq:poly_drift}
  \norm{ \dvar(t+1) - \dvar^\ast } - \norm{ \dvar(t) - \dvar^\ast } \leq - \frac{L_\poly}{2V}.
\end{equation}  
\end{lemma}
\begin{IEEEproof}
From Lemma \ref{lem:dual_move}, suppose the following condition holds
\begin{equation}
  \label{eq:p_1}
   \frac{2}{V} [d(\dvar(t)) - d(\dvar^\ast)] + \frac{2 C}{V^2} \leq -\frac{L_\poly}{V} \norm{ \dvar(t) - \dvar^\ast } + \frac{ L_\poly^2 }{4 V^2},
\end{equation}
then inequality \eqref{eq:dual_move} becomes
\begin{align*}
  \norm{\dvar(t+1) - \dvar^\ast}^2
  & \leq \norm{ \dvar(t) - \dvar^\ast }^2 - \frac{L_\poly}{V} \norm{ \dvar(t) - \dvar^\ast } + \frac{ L_\poly^2 }{4 V^2} \\
  & = \prts{ \norm{ \dvar(t) - \dvar^\ast } - \frac{L_\poly}{2V} }^2.
\end{align*}
It follows that if $\norm{ \dvar(t) - \dvar^\ast } \geq B_\poly(V) \geq \frac{L_\poly}{2V}$, then inequality \eqref{eq:poly_drift} holds.

It requires to show that condition \eqref{eq:p_1} holds when $\norm{ \dvar(t) - \dvar^\ast } \geq B_\poly(V)$.  Note that condition \eqref{eq:p_1} holds when
\begin{equation*}
  d(\dvar(t)) - d(\dvar^\ast) \leq - \frac{C}{V} - \frac{L_\poly}{2} \norm{ \dvar(t) - \dvar^\ast }.
\end{equation*}
By the locally-polyhedral property \eqref{eq:poly}, if $-L_\poly \norm{ \dvar(t) - \dvar^\ast } \leq - \frac{C}{V} - \frac{L_\poly}{2} \norm{ \dvar(t) - \dvar^\ast }$, then the above inequality holds.  This means that condition \eqref{eq:p_1} holds when $\norm{\dvar(t) - \dvar^\ast} \geq \frac{2 C}{V L_\poly}$.  This proves the lemma.
\end{IEEEproof}

Lemma \ref{lem:poly_drift} implies that, if the distance between $\dvar(t)$ and $\dvar^\ast$ is at least $B_\poly(V)$, the successor $\dvar(t+1)$ will be closer to $\dvar^\ast$.  This suggests the existence of a convergence set in which a subsequence of $\{ \dvar(t) \}_{t=0}^\infty$ resides.  Note that $\sqrt{2C}/V$ bounds $\norm{ \dvar(t+1) - \dvar(t) }$ for all $t$ as in \eqref{eq:onestep_bound}.  

The steady state of Algorithm \ref{alg:Subgrad} is defined from this set.  This convergence set is defined as
\begin{equation}
  \label{eq:poly_region}
  \set{R}_\poly(V) = \prtc{ \dvar \in \Pi : \norm{ \dvar - \dvar^\ast } \leq B_\poly(V) + \frac{\sqrt{2 C}}{V} }.
\end{equation}

Let $T_\poly$ be the first iteration that a generated dual variable enters this set:
\begin{equation}
  \label{eq:poly_time}
  T_\poly = \arginf_{t \geq 0} \prtc{ \dvar(t) \in \set{R}_\poly(V) }.
\end{equation}
Intuitively, $T_\poly$ is the end of the transient phase and is the beginning of the steady state phase.

\begin{lemma}
  Under Assumptions \ref{ass:uniqueness} and \ref{ass:poly}, $T_\poly \leq O(V)$.
\end{lemma}
\begin{IEEEproof}
  Since $\norm{ \dvar(0) - \dvar^\ast }$ is a constant, Lemma \ref{lem:poly_drift} proves the claim.
\end{IEEEproof}

Then we show that dual variables generated after iteration $T_\poly$ never leave $\set{R}_\poly(V)$.
\begin{lemma}
  \label{lem:poly_region}
  Under Assumptions \ref{ass:uniqueness} and \ref{ass:poly}, the generated dual variables from Algorithm \ref{alg:Subgrad} satisfy $\dvar(t) \in \set{R}_\poly(V)$ for all $t \geq T_\poly$.
\end{lemma}
\begin{IEEEproof}
We prove the lemma by induction.  
First we note that $\dvar(T_\poly) \in \set{R}_\poly(V)$ by the definition of $T_\poly$.
Suppose that $\dvar(t) \in \set{R}_\poly(V)$.
Then two cases are considered. 

i) If $\norm{\dvar(t) - \dvar^\ast} \geq B_\poly(V)$, it follows from \eqref{eq:poly_drift} that
\begin{equation*}
\norm{\dvar(t+1) - \dvar^\ast} \leq \norm{\dvar(t) - \dvar^\ast} - \frac{L_\poly}{2V} \leq B_\poly(V) + \frac{\sqrt{2 C}}{V}.
\end{equation*}

ii) If $\norm{\dvar(t) - \dvar^\ast} \leq B_\poly(V)$, it follows from the triangle inequality that 
\begin{align*}
  \norm{\dvar(t+1) - \dvar^\ast} & \leq \norm{\dvar(t+1) - \dvar(t)} + \norm{\dvar(t) - \dvar^\ast} \\
  & \leq \frac{\sqrt{2 C}}{V} + B_\poly(V),
\end{align*}
by \eqref{eq:onestep_bound} and the assumption of $\norm{ \dvar(t) - \dvar^\ast }$.
Hence, $\dvar(t+1) \in \set{R}_\poly(V)$ in both cases.
This proves the lemma by induction.
\end{IEEEproof}

Finally, a convergence result is ready to be stated.
Let $\overline{a_{T_\poly}}(T) = \frac{1}{T} \sum_{t=T_\poly}^{T_\poly+T-1} a(t)$ be an average of sequence $\{a(t)\}_{t=T_\poly}^{T_\poly+T-1}$ that starts from $T_\poly$.
\begin{theorem}
  \label{thm:poly}
  Under Assumptions \ref{ass:uniqueness} and \ref{ass:poly}, for $T > 0$, let $\prtc{ x(t), w(t) }_{t=T_\poly}^\infty$ be a subsequence generated by Algorithm \ref{alg:Subgrad}, where $T_\poly$ is defined in \eqref{eq:poly_time}.  The following bounds hold:
\begin{align}
  f(\overline{x_{T_\poly}}(T)) - f^\optvar
  & \leq \frac{C}{V} + \frac{2 V M}{T} \prtbs{ \frac{\sqrt{2C}}{V} + B_\poly(V) } \notag\\
  & \quad + \frac{V}{2T} \biggl\{ \prtbs{ \frac{\sqrt{2C}}{V} + B_\poly(V) }^2 + 4 \norm{\dvar^\ast}\prts{ \frac{\sqrt{2C}}{V} + B_\poly(V) } \biggr\}   \label{eq:poly_avg}\\
  g_j(\overline{x_{T_\poly}}(T)) & \leq \frac{2 V(1+M)}{T} \prtbs{ \frac{\sqrt{2C}}{V} + B_\poly(V) }, \quad j \in \{ 1, \dotsc, J \}.\label{eq:poly_vio}
\end{align}
\end{theorem}
\begin{IEEEproof}
The first part of the theorem follows from \eqref{eq:avg_bound} with the average starting from $T_\poly$ that
\begin{equation}
  \label{eq:pt_1}
  f(\overline{x_{T_\poly}}(T)) - f^\optvar \leq \frac{C}{V} +  \frac{V}{2T} \prts{ \norm{\dvar(T_\poly)}^2 - \norm{\dvar(T_\poly + T)}^2 } \\+ \frac{V M}{T} \norm{ z(T_\poly + T) - z(T_\poly) }.
\end{equation}

For any $\dvar \in \Pi$, it holds that: 
\begin{equation*}
\norm{\dvar}^2 = \norm{ \dvar - \dvar^\ast }^2 + \norm{\dvar^\ast}^2 + 2[\dvar - \dvar^\ast]^\tr \dvar^\ast.
\end{equation*}
The second term on the right-hand-side of \eqref{eq:pt_1} can be upper bounded by applying this equality.
\begin{align}
  \norm{\dvar(T_\poly)}^2 - \norm{\dvar(T_\poly + T)}^2
  & = \norm{ \dvar(T_\poly) - \dvar^\ast }^2 + 2[\dvar(T_\poly) - \dvar^\ast]^\tr \dvar^\ast \notag\\
  & \quad - \norm{ \dvar(T_\poly + T) - \dvar^\ast }^2 - 2[\dvar(T_\poly + T) - \dvar^\ast]^\tr \dvar^\ast \notag\\
  & \leq \norm{\dvar(T_\poly) - \dvar^\ast}^2 + 2 [\dvar(T_\poly) - \dvar(T_\poly+T)]^\tr \dvar^\ast \notag\\
  & \leq \norm{\dvar(T_\poly) - \dvar^\ast}^2 + 2 \norm{\dvar(T_\poly) - \dvar(T_\poly+T)} \norm{\dvar^\ast} \label{eq:pt_2}
\end{align}
From Lemma \ref{lem:poly_region}, the first term of \eqref{eq:pt_2} is bounded by $\norm{\dvar(T_\poly) - \dvar^\ast}^2 \leq [\sqrt{2C}/V + B_\poly(V)]^2$.  
From triangle inequality and Lemma \ref{lem:poly_region}, the last term of \eqref{eq:pt_2} is bounded by
\begin{align}
  \norm{\dvar(T_\poly+T) - \dvar(T_\poly)}
  & \leq \norm{\dvar(T_\poly+T) - \dvar^\ast} + \norm{\dvar^\ast - \dvar(T_\poly)} \notag\\
  & \leq 2 \prtbs{ \sqrt{2C}/V + B_\poly(V) }. \label{eq:pt_3}
\end{align}
Therefore, inequality \eqref{eq:pt_2} is bounded from above by $[\sqrt{2C}/V + B_\poly(V)]^2 + 4 \norm{\dvar^\ast} [\sqrt{2C}/V + B_\poly(V)]$.
Substituting this bound into \eqref{eq:pt_1} and using the fact that
\begin{equation*}
  \norm{ z(T_\poly+T) - z(T_\poly) } \leq \norm{ \dvar(T_\poly + T) - \dvar(T_\poly) } \leq 2[ \sqrt{2C}/V + B_\poly(V) ]
\end{equation*}
proves the first part of the theorem.

The last part follows from \eqref{eq:avg_vio} that
\begin{equation*}
  g_j(\overline{x_{T_\poly}}(T)) \leq \frac{V}{T} \abs{ w_j(T_\poly+T) - w_j(T_\poly) } \\+ \frac{V M}{T} \norm{ z(T_\poly+T) - z(T_\poly) }.
\end{equation*}
Since $\abs{ w_j(T_\poly+T) - w_j(T_\poly) }$ and $\norm{ z(T_\poly+T) - z(T_\poly) }$ are bounded above by $\norm{\dvar(T_\poly+T) - \dvar(T_\poly)}$, the above inequality is upper bounded by
\begin{align*}
  g_j(\overline{x_{T_\poly}}(T)) & \leq \frac{V(1+M)}{T} \norm{\dvar(T_\poly+T) - \dvar(T_\poly)} \\
  & \leq \frac{2 V(1+M)}{T} \prtbs{ \frac{\sqrt{2C}}{V} + B_\poly(V) },
\end{align*}
where the last inequality uses relation \eqref{eq:pt_3}.  
This proves the last part of the theorem.
\end{IEEEproof}

Theorem \ref{thm:poly} can be interpreted as follows.  The deviation from the optimality value \eqref{eq:poly_avg} is bounded above by $O(1/V + 1/T)$.  The constraint violation \eqref{eq:poly_vio} is bounded above by $O(1/T)$.  To have both bounds be within $O(\epsilon)$, we set $V = 1/\epsilon$ and $T = 1/\epsilon$, and the convergence time of Algorithm \ref{alg:Subgrad} is $O(1/\epsilon)$.  Note that both bounds consider the average starting after reaching the steady state at time $T_\poly$, and this transient time $T_\poly$ is at most $O(1/\epsilon)$.

\subsection{Locally-Smooth Dual Function}
\label{sec:smooth}

Throughout Section \ref{sec:smooth}, the dual function \eqref{eq:dual_function} is assumed to have a locally-smooth property, introduced in \cite{Longbo:Delay_reduction}, as stated in Assumption \ref{ass:smooth} and illustrated in Figure \ref{fig:ps_functions}.
\begin{assumption}
  \label{ass:smooth}
  Let $\dvar^\ast$ be the unique Largrange multiplier, there exist $S > 0$ and $L_\smooth > 0$ such that whenever $\dvar \in \Pi$ and $\norm{ \dvar - \dvar^\ast } \leq S$, dual function \eqref{eq:dual_function} satisfies
\begin{equation}
  \label{eq:smooth}
  d(\dvar^\ast) \geq d(\dvar) + L_\smooth \norm{ \dvar - \dvar^\ast }^2.
\end{equation}
Also, there exists $D_\smooth > 0$ such that whenever $\dvar \in \Pi$ and $d(\dvar^\ast) - d(\dvar) \leq D_\smooth$, dual variable satisfies $\norm{ \dvar - \dvar^\ast } \leq S$.
\end{assumption}
The ``s'' subscript in $L_\smooth$ represents ``smooth.''

The behavior of the generated dual variables from a dual function satisfying the locally-smooth assumption can be described as follows. Define 
\begin{equation*}
B_\smooth(V) \defeq \max\prtc{ \frac{1}{V^{1.5}}, \frac{\sqrt{V} + \sqrt{V + 4 L_\smooth C V}}{2 L_\smooth V} }.
\end{equation*}

\begin{lemma}
  \label{lem:smooth_drift}
  Under Assumptions \ref{ass:uniqueness} and \ref{ass:smooth}, for sufficiently large $V$ that $B_\smooth(V) < S$, whenever $B_\smooth(V) \leq \norm{ \dvar(t) - \dvar^\ast } \leq S$, it follows that
\begin{equation}
  \label{eq:smooth_drift}
  \norm{ \dvar(t+1) - \dvar^\ast } - \norm{ \dvar(t) - \dvar^\ast } \leq - \frac{1}{V^{1.5}}.
\end{equation}
\end{lemma}

\begin{IEEEproof}
From Lemma \ref{lem:dual_move}, suppose the following condition holds
\begin{equation}
  \label{eq:s_2}
  \frac{2}{V} [d(\dvar(t)) - d(\dvar^\ast)] + \frac{2 C}{V^2} \leq -\frac{2}{V^{1.5}} \norm{ \dvar(t) - \dvar^\ast } + \frac{1}{V^3},
\end{equation}
then inequality \eqref{eq:dual_move} becomes
\begin{align*}
  \norm{ \dvar(t+1) - \dvar^\ast }^2 
  & \hspace{-1mm}\leq \hspace{-1mm} \norm{ \dvar(t) - \dvar^\ast }^2 - \frac{2}{V^{1.5}} \norm{ \dvar(t) - \dvar^\ast } + \frac{1}{V^3} \\
  & \hspace{-1mm}= \hspace{-1mm} \prts{ \norm{ \dvar(t) - \dvar^\ast } - \frac{1}{V^{1.5}} }^2.
\end{align*}
Furthermore, if $\norm{ \dvar(t) - \dvar^\ast } \geq B_\smooth(V) \geq \frac{1}{V^{1.5}}$, then the desired inequality \eqref{eq:smooth_drift} holds.

It requires to show that condition \eqref{eq:s_2} holds when $S \geq \norm{ \dvar(t) - \dvar^\ast } \geq B_\smooth(V)$.  Condition \eqref{eq:s_2} holds when 
\begin{equation*}
d(\dvar(t)) - d(\dvar^\ast) \leq - \frac{C}{V} - \frac{1}{\sqrt{V}} \norm{ \dvar(t) - \dvar^\ast }.
\end{equation*}
By the locally-smooth property \eqref{eq:smooth}, if $ -L_\smooth \norm{ \dvar(t) - \dvar^\ast }^2 \leq - \frac{C}{V} - \frac{1}{\sqrt{V}} \norm{ \dvar(t) - \dvar^\ast }$, then the above inequality holds.  
This means that condition \eqref{eq:s_2} holds when
\begin{equation*}
  L_\smooth \norm{ \dvar(t) - \dvar^\ast }^2 - \frac{1}{\sqrt{V}} \norm{ \dvar(t) - \dvar^\ast } - \frac{C}{V} \geq 0.
\end{equation*}
The above inequality happens when 
\begin{align*}
  \norm{ \dvar(t) - \dvar^\ast } & \geq \frac{ \frac{1}{\sqrt{V}} + \sqrt{ \frac{1}{V} + 4 L_\smooth \frac{C}{V} }}{2L_\smooth} \\
  & = \frac{\sqrt{V} + \sqrt{V + 4 L_\smooth C V}}{2 L_\smooth V}.
\end{align*}
This prove the lemma.
\end{IEEEproof}

Lemma \ref{lem:smooth_drift} suggests the existence of a convergence set.
The steady state of Algorithm \ref{alg:Subgrad} is also defined from this set as
\begin{equation}
  \label{eq:smooth_region}
  \set{R}_\smooth(V) = \prtc{ \dvar \in \Pi : \norm{ \dvar - \dvar^\ast } \leq B_\smooth(V) + \frac{\sqrt{2 C}}{V} }.
\end{equation}
Let $T_\smooth$ denote the first iteration that a generated dual variables arrives at the convergence set:
\begin{equation}
  \label{eq:smooth_time}
T_\smooth = \arginf_{t \geq 0} \prtc{ \dvar(t) \in \set{R}_\smooth(V) }.
\end{equation}

\begin{lemma}
  Under Assumptions \ref{ass:uniqueness} and \ref{ass:smooth}, when $V$ is sufficiently large and $B_\smooth(V) < S$, then $T_\smooth \leq O(V^{1.5})$.
\end{lemma}
\begin{IEEEproof}
We first shows that there exists $t' \leq O(V)$ such that $\norm{\dvar(t') - \dvar^\ast} \leq S$.  We show that the following is true:
\begin{equation}
  \label{eq:transient_dual}
  d(\dvar^\ast) - \max_{0 \leq t \leq E_\delta(V)} d(\dvar(t)) \leq \frac{C}{V} + \frac{\delta}{2},
\end{equation}
where $E_\delta(V) \defeq \left\lfloor \frac{V \norm{\dvar(0) - \dvar^\ast}^2}{\delta} \right\rfloor$.

This is proved by contradiction.  Suppose inequality \eqref{eq:transient_dual} does not hold, i.e., 
\begin{equation*}
d(\dvar^\ast) - d(\dvar(t)) > \frac{C}{V} + \frac{\delta}{2} \quad \text{for all}~ 0 \leq t \leq E_\delta(V).
\end{equation*}
From inequality \eqref{eq:dual_move}, it follows that for $0 \leq t \leq E_\delta(V)$
\begin{align*}
  \norm{\dvar(t+1) - \dvar^\ast}^2 
  & \leq \norm{ \dvar(t) - \dvar^\ast }^2 + \frac{2 C}{V^2} - \frac{2}{V} \prtr{ \frac{C}{V} + \frac{\delta}{2} } \\
  & \leq \norm{ \dvar(t) - \dvar^\ast }^2 - \frac{\delta}{V}.
\end{align*}
Summing from $t = 0, \dotsc, E_\delta(V)$ yields:
\begin{equation*}
  \norm{\dvar(E_\delta(V)+1) - \dvar^\ast}^2 \leq \norm{\dvar(0) - \dvar^\ast}^2 - \frac{ [E_\delta(V)+1] \delta}{V},
\end{equation*}
and $E_\delta(V)+1 \leq \frac{V \norm{\dvar(0) - \dvar^\ast}^2}{\delta}$.  This contradicts the definition of $E_\delta(V)$.  Thus, property \eqref{eq:transient_dual} holds.

Let $\delta = D_\smooth$ and $V > 2 C /D_\smooth$, we have $d(\dvar^\ast) - d(\dvar(t)) \leq D_\smooth$ for some $0 \leq t \leq E_\delta(V)$.  
Then from Assumption \ref{ass:smooth}, we have $\norm{\dvar(t) - \dvar^\ast} \leq S$, and by the definition of $E_\delta(V)$, it takes at most $O(V)$ to arrive where the locally-smooth assumption holds.  Then Lemma \ref{lem:smooth_drift} implies that the algorithm needs at most $O(V^{1.5})$ to enter the convergence set.
\end{IEEEproof}

Next we show that, once the sequence of dual variables enters $\set{R}_\smooth(V)$, it never leaves the set.

\begin{lemma}
  \label{lem:smooth_region}
  Under Assumptions \ref{ass:uniqueness} and \ref{ass:smooth}, when $V$ is sufficiently large and $B_\smooth(V) + \frac{\sqrt{2 C}}{V} < S$, the generated dual variables from Algorithm \ref{alg:Subgrad} satisfy $\dvar(t) \in \set{R}_\smooth(V)$ for all $t \geq T_\smooth$.
\end{lemma}
\begin{IEEEproof}
We prove the lemma by induction.  
First we note that $\dvar(T_\smooth) \in \set{R}_\smooth(V)$ by its definition.  
Suppose that $\dvar(t) \in \set{R}_\smooth(V)$, which implies that $\norm{ \dvar(t) - \dvar^\ast } \leq B_\smooth(V) + \sqrt{2C}/V < S$.  
Then two cases are considered. 

i) If $\norm{\dvar(t) - \dvar^\ast} > B_\smooth(V)$, it follows from \eqref{eq:smooth_drift} that
\begin{equation*}
  \norm{\dvar(t+1) - \dvar^\ast} < \norm{\dvar(t) - \dvar^\ast} - \frac{1}{V^{1.5}} < B_\smooth(V) + \frac{\sqrt{2 C}}{V}.  
\end{equation*}

ii) If $\norm{\dvar(t) - \dvar^\ast} \leq B_\smooth(V)$, it follows from the triangle inequality and \eqref{eq:onestep_bound} that 
\begin{align*}
  \norm{\dvar(t+1) - \dvar^\ast} & \leq \norm{\dvar(t+1) - \dvar(t)} + \norm{\dvar(t) - \dvar^\ast} \\
  & \leq \frac{\sqrt{2 C}}{V} + B_\smooth(V).  
\end{align*}
Hence, $\dvar(t+1) \in \set{R}_\smooth(V)$ in both cases.  
This proves the lemma by induction.
\end{IEEEproof}

Now a convergence of a steady state is ready to be stated.  
\begin{theorem}
  \label{thm:smooth}
  Under Assumptions \ref{ass:uniqueness} and \ref{ass:smooth}, when $V$ is sufficiently large and $B_\smooth(V) + \frac{\sqrt{2 C}}{V} < S$, for $T > 0$, let $\prtc{ x(t), w(t) }_{t=T_\smooth}^\infty$ be a subsequence generated by Algorithm \ref{alg:Subgrad}, where $T_\smooth$ is defined in \eqref{eq:smooth_time}.  
The following bounds hold:
\begin{align}
  f(\overline{x_{T_\smooth}}(T)) - f^\optvar & \leq \frac{C}{V} + \frac{2 V M}{T} \prtbs{ \frac{\sqrt{2C}}{V} + B_\smooth(V) }\notag\\
  & \quad + \frac{V}{2T} \biggl\{ \prtbs{ \frac{\sqrt{2C}}{V} + B_\smooth(V) }^2 + 2 \norm{\dvar^\ast}\prts{ \frac{\sqrt{2C}}{V} + B_\smooth(V) } \biggr\}   \label{eq:smooth_avg} \\
  g_j(\overline{x_{T_\smooth}}(T)) & \leq \frac{2 V(1+M)}{T} \prtbs{ \frac{\sqrt{2C}}{V} + B_\smooth(V) }, \quad j \in \{ 1, \dotsc, J \}.\label{eq:smooth_vio}
\end{align}
\end{theorem}
\begin{IEEEproof}
The first part of the theorem follows from \eqref{eq:avg_bound} with the average starting from $T_\smooth$ that
\begin{equation}
  \label{eq:st_1}
  f(\overline{x_{T_\smooth}}(T)) - f^\optvar \leq \frac{C}{V} +  \frac{V}{2T} \prts{ \norm{\dvar(T_\smooth)}^2 - \norm{\dvar(T_\smooth + T)}^2 } \\+ \frac{V M}{T} \norm{ z(T_\smooth + T) - z(T_\smooth) }.
\end{equation}

The second term on the right-hand-side of \eqref{eq:st_1} can be bounded from above by
\begin{align}
  \norm{\dvar(T_\smooth)}^2 - \norm{\dvar(T_\smooth + T)}^2
  &  \leq \prtbs{ \frac{\sqrt{2C}}{V} + B_\smooth(V) }^2 + 4\norm{\dvar^\ast} \prtbs{ \frac{\sqrt{2C}}{V} + B_\smooth(V) }, \label{eq:st_2}
\end{align}
where the deviation is similar to steps in \eqref{eq:pt_2} and \eqref{eq:pt_3}.

The last term on the right-hand-side of \eqref{eq:st_1} can be bounded from above by
\begin{align}
  \norm{z(T_\smooth+T) - z(T_\smooth)} & \leq 2 \prtbs{ \sqrt{2C}/V + B_\smooth(V) }. \label{eq:st_3}
\end{align}
Substituting bounds \eqref{eq:st_2} and \eqref{eq:st_3} into \eqref{eq:st_1} proves the first part of the theorem.

The last part follows from \eqref{eq:avg_vio} that
\begin{equation*}
  g_j(\overline{x_{T_\smooth}}(T)) \leq \frac{V}{T} \abs{ w_j(T_\smooth+T) - w_j(T_\smooth) } \\+ \frac{V M}{T} \norm{ z(T_\smooth+T) - z(T_\smooth) }.
\end{equation*}
Since $\abs{ w_j(T_\smooth+T) - w_j(T_\smooth) }$ and $\norm{ z(T_\smooth+T) - z(T_\smooth) }$ are bounded above by $\norm{\dvar(T_\smooth+T) - \dvar(T_\smooth)}$, the above inequality is upper bounded by
\begin{align*}
  g_j(\overline{x_{T_\smooth}}(T)) & \leq \frac{V(1+M)}{T} \norm{\dvar(T_\smooth+T) - \dvar(T_\smooth)} \\
  & \leq \frac{2 V(1+M)}{T} \prtbs{ \frac{\sqrt{2C}}{V} + B_\smooth(V) }.
\end{align*}
This proves the last part of the theorem.
\end{IEEEproof}

Theorem \ref{thm:smooth} can be interpreted as follows.  
The deviation from the optimality \eqref{eq:smooth_avg} is bounded above by $O(1/V + \sqrt{V}/T)$.  The constraint violation \eqref{eq:smooth_vio} is bounded above by $O(\sqrt{V}/T)$.  To have both bounds be within $O(\epsilon)$, we set $V = 1/\epsilon$ and $T = 1/{\epsilon^{1.5}}$, and the convergence time of Algorithm \ref{alg:Subgrad} is $O(1/{\epsilon^{1.5}})$.  Note that both bounds consider the average starting after reaching the steady state at time $T_\smooth$, and this transient time $T_\smooth$ is at most $O(1/{\epsilon^{1.5}})$.

\subsection{Staggered Time Averages}

In order to take advantage of the improved convergence rates, computing time averages must be started after the transient phase.  To achieve this performance without determining the exact end time of the transient phase, time averages can be restarted over successive frames whose frame lengths increase geometrically.  For example, if one triggers a restart at times $2^k$ for integers $k$, then a restart is guaranteed to occur within a factor of $2$ of the time of the actual end of the transient phase.

\subsection{Summary of Convergence Results}
The results in Theorems \ref{thm:gen_conv}, \ref{thm:poly}, and \ref{thm:smooth} (denoted by General, Polyhedron, and Smooth) are summarized in Table \ref{tlb:convergence}.  
Note that the general convergence time is considered to be in the steady state from the beginning.

\begin{table}
  \centering
  \caption{Convergence Times}
  \label{tlb:convergence}
  \begin{tabular}{||l||c|c|c||}
    \hline\hline
    & General & Polyhedron & Smooth \\
    \hline\hline
    Transient state & $0$ & $O(1/{\epsilon})$ & $O(1/{\epsilon^{1.5}})$ \\
    \hline
    Steady state & $O(1/\epsilon^2)$ & $O(1/{\epsilon})$ & $O(1/{\epsilon^{1.5}})$ \\
    \hline\hline
  \end{tabular}
\end{table}

\section{Sample Problems}
\label{sec:example}

This section illustrates the convergence times of the time-average Algorithm \ref{alg:Subgrad} under locally-polyhedral and locally-smooth assumptions.  A considered formulation is
\begin{align}
  \label{eq:sim_formulation}
  \minimize \quad & f(\bar{x}) \\
  \subjectto \quad & 2 \bar{x}_1 + \bar{x}_2 \geq 1.5, \quad \bar{x}_1 + 2 \bar{x}_2 \geq 1.5 \notag\\
  & x_1(t), x_2(t) \in \{0, 1, 2, 3\}, \quad t \in \{0, 1, 2, \dotsc\} \notag
\end{align}
where function $f$ will be given for different cases.

\begin{figure}
  \centering
  \includegraphics[scale=0.43]{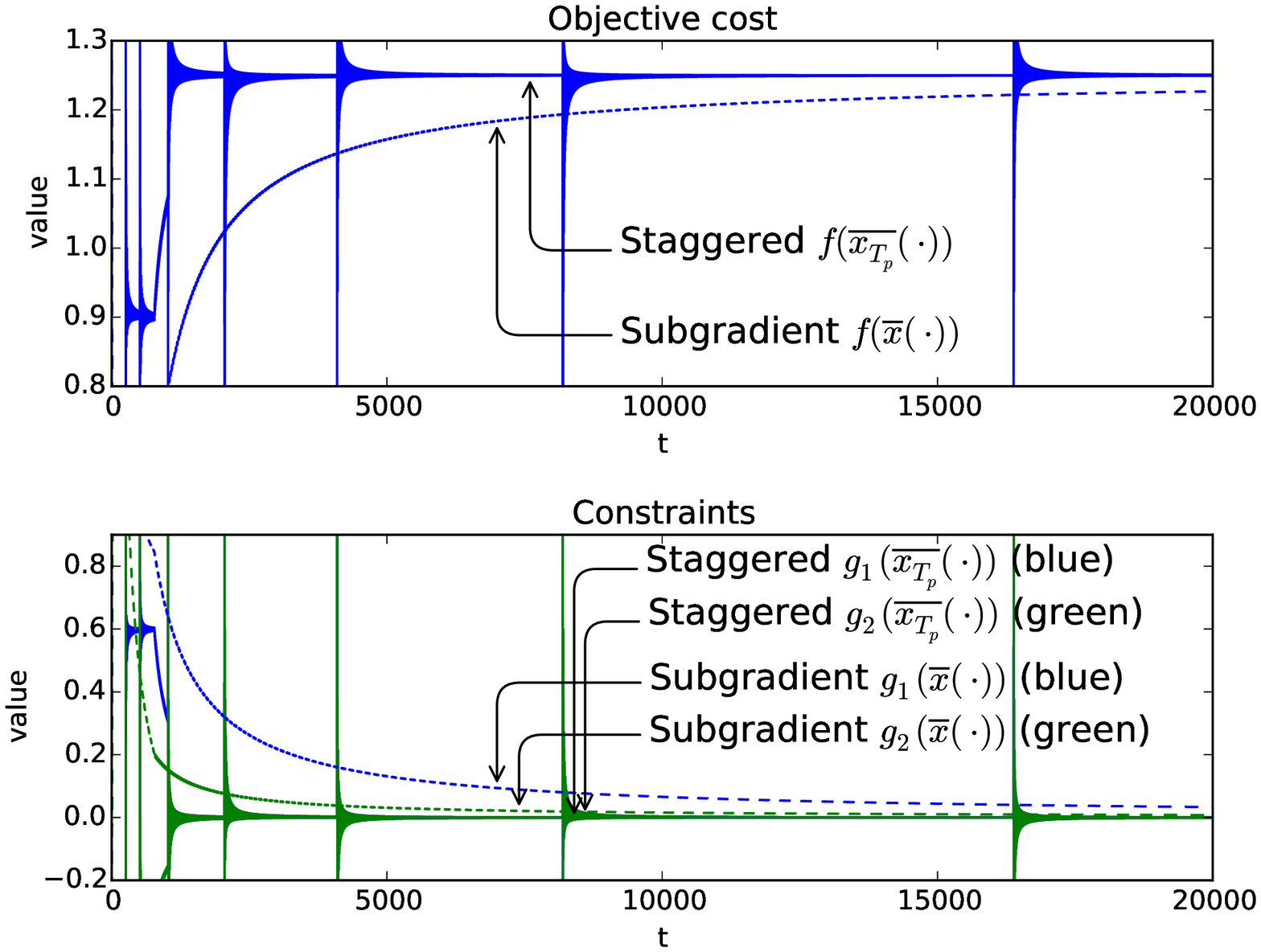}
  \caption{Iterations solving problem \eqref{eq:sim_formulation} with $f(x) = 1.5 x_1 + x_2$}
  \label{fig:poly_w_ass}
\end{figure}

Under the locally-polyhedral assumption, let $f(x) = 1.5 x_1 + x_2$ be the objective function of problem \eqref{eq:sim_formulation}.  In this setting, the optimal value is $1.25$ when $\bar{x}_1 = \bar{x}_2 = 0.5$.  Figure \ref{fig:poly_w_ass} shows the values of objective and constraint functions of time-averaged solutions.  It is easy to see the faster convergence time $O(1/\epsilon)$ from the polyhedral result ($T_\poly = 2048$) compared to a general result with convergence time $O(1/\epsilon^2)$.

\begin{figure}
  \centering
  \includegraphics[scale=0.43]{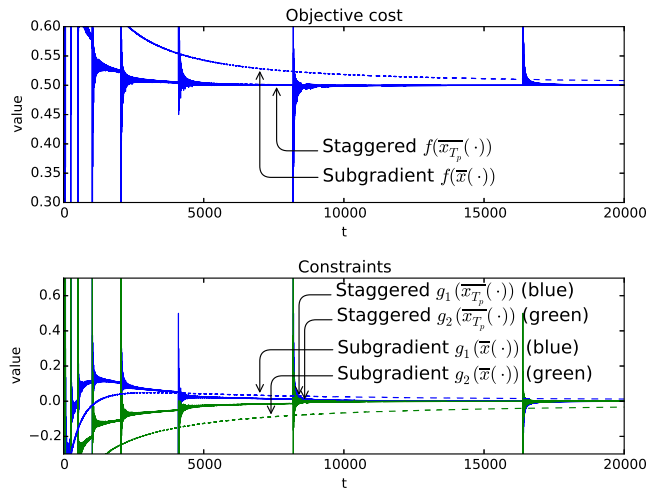} 
  \caption{Iterations solving problem \eqref{eq:sim_formulation} with $f(x) = x_1^2 + x_2^2$}
  \label{fig:smooth_w_ass}
\end{figure}

Under the locally-smooth assumption, let $f(x) = x_1^2 + x_2^2$ be the objective function of problem \eqref{eq:sim_formulation}.  Note that the optimal value of this problem is $0.5$ where $\bar{x}_1 = \bar{x}_2 = 0.5$.  Figure \ref{fig:smooth_w_ass} shows the values of objective and constraint functions of time-averaged solutions.  The smooth result starts the average from $(T_\smooth=)8192^\nth$ iterations.  It is easy to see that the general result converges slower than the smooth result.  This illustrates the difference between $O(1/\epsilon^2)$ and $O(1/{\epsilon^{1.5}})$.

\begin{figure}
  \centering
  \includegraphics[scale=0.43]{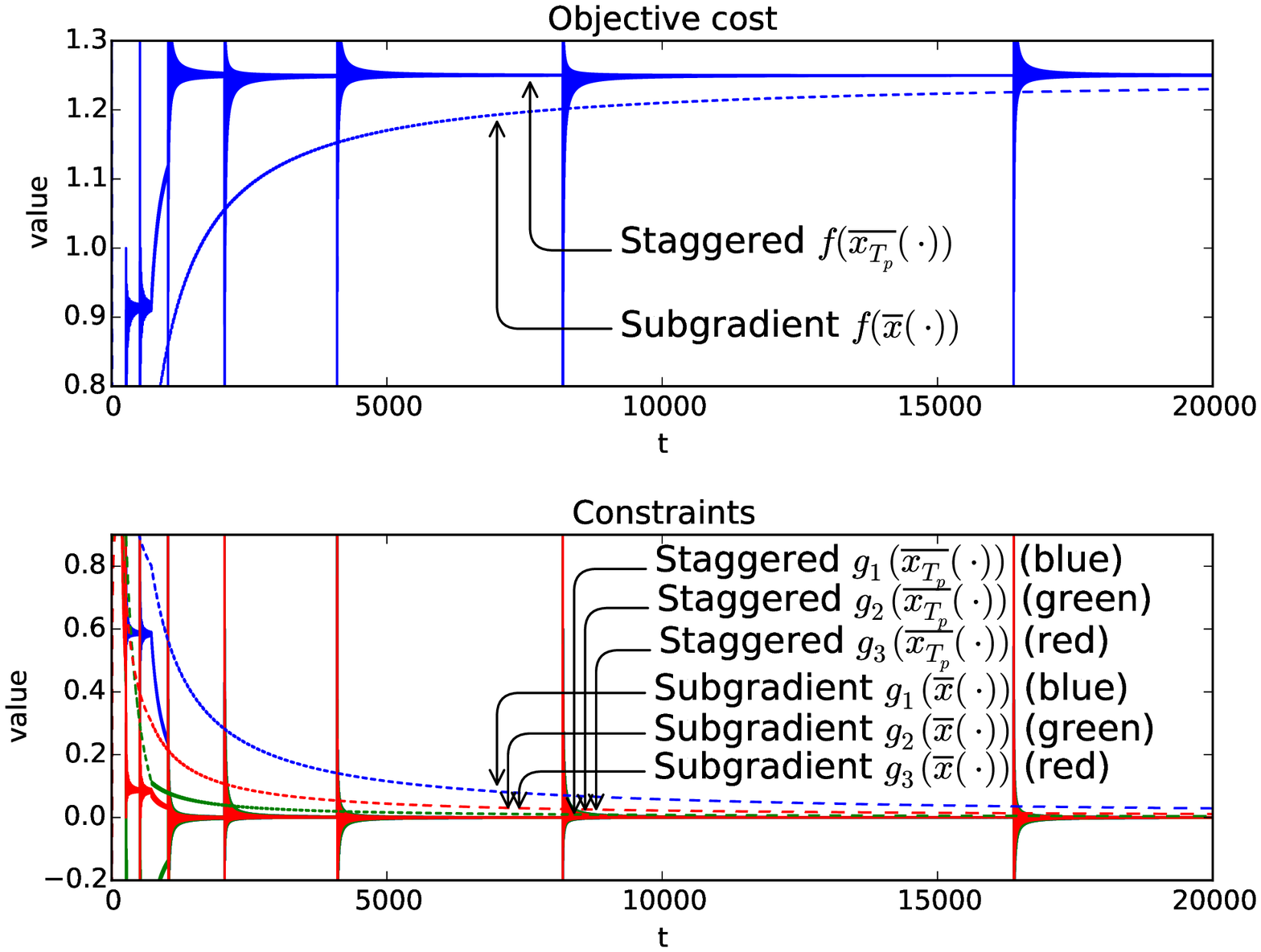} 
  \caption{Iterations solving problem \eqref{eq:sim_formulation} with $f(x) = 1.5 x_1 + x_2$ and an additional constraint $\bar{x}_1 + \bar{x}_2 \geq 1$}
  \label{fig:poly_wo_ass}
\end{figure}

\begin{figure}
  \centering
  \includegraphics[scale=0.43]{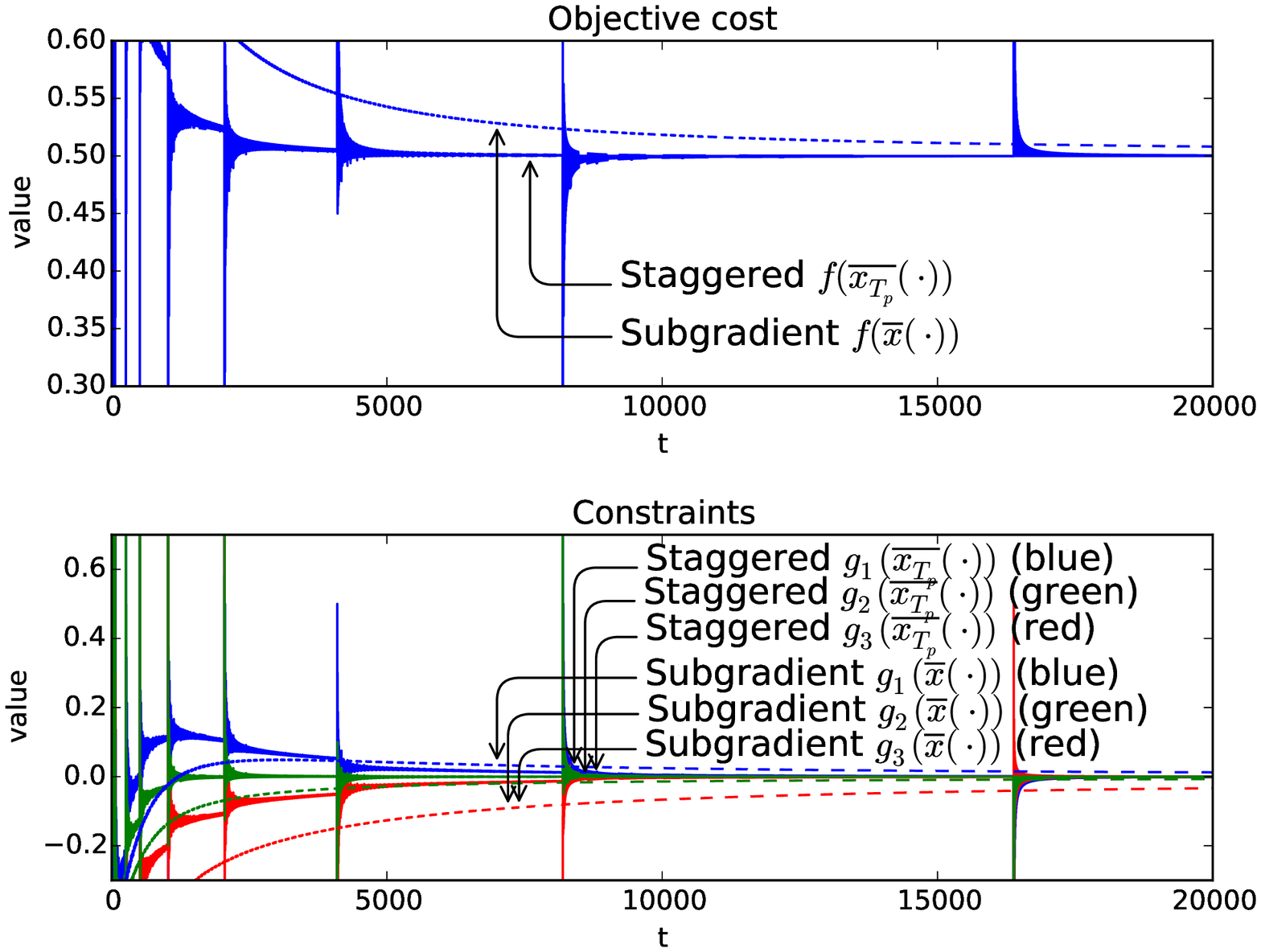} 
  \caption{Iterations solving problem \eqref{eq:sim_formulation} with $f(x) = x_1^2 + x_2^2$ and an additional constraint $\bar{x}_1 + \bar{x}_2 \geq 1$}
  \label{fig:smooth_wo_ass}
\end{figure}

Figures \ref{fig:poly_wo_ass} and \ref{fig:smooth_wo_ass} illustrate the convergence times of problems, defined in each figure's caption, without the uniqueness assumption.  The Comparison of Figures \ref{fig:poly_w_ass} and \ref{fig:poly_wo_ass} shows that there is no difference in the order of convergence time.  Similarly, figures \ref{fig:smooth_w_ass} and \ref{fig:smooth_wo_ass} show no difference in terms of the order of convergence.

\section{CONCLUSION}
\label{sec:conclusion}
We consider the time-average optimization problem with a non-convex (possibly discrete) decision set.  We show that the problem has a corresponding (one-shot) convex optimization formulation.  This connects the Lyapunov optimization technique and convex optimization theory.  Using convex analysis we prove a general convergence time of $O(1/\epsilon^2)$ when the Slater condition holds.  Under an assumption on the uniqueness of a Lagrange multiplier, we prove that faster convergence times $O(1/{\epsilon})$ and $O(1/{\epsilon^{1.5}})$ are possible for locally-polyhedral and locally-smooth problems.

\bibliographystyle{IEEEtran}
\bibliography{Reference}

\end{document}